\documentclass{ert-like}
\usepackage{enumerate, amsfonts, amssymb}

\newlength{\tabwidth}
\newlength{\tabheight}
\newlength{\tabrule}
\newlength{\tabwidthx}
\newlength{\tabheightx}

\def\gentabbox#1#2#3#4{\vbox to \tabheight{\setlength{\tabrule}{#3}%
  \setlength{\tabwidthx}{#1\tabwidth}\addtolength{\tabwidthx}{\tabrule}%

\setlength{\tabheightx}{#2\tabheight}\addtolength{\tabheightx}{-\tabheight}%
  \hbox to #1\tabwidth{%
 \hspace{-0.5\tabrule}\rule{\tabrule}{#2\tabheight}\hspace{-\tabrule}%
    \vbox to #2\tabheight{\hsize=\tabwidthx%
      \vspace{-0.5\tabrule}\hrule width\tabwidthx height\tabrule%
      \vspace{-0.5\tabrule}\vfil%
      \hbox to \tabwidthx{\hss#4\hss}%
        \vfil\vspace{-0.5\tabrule}%
      \hrule width\tabwidthx height\tabrule\vspace{-0.5\tabrule}}%
 \hspace{-\tabrule}\rule{\tabrule}{#2\tabheight}\hspace{-0.5\tabrule}}%
  \vspace{-\tabheightx}}}
\def\genblankbox#1#2{\vbox to \tabheight{\vfil\hbox to
#1\tabwidth{\hfil}}}
\def\tabbox#1#2#3{\gentabbox{#1}{#2}{0.4pt}{\strut #3}}

\catcode`\:=13 \catcode`\.=13 \catcode`\;=13 
\catcode`\>=13 \catcode`\^=13
\def:#1\\{\hbox{$#1$}}
\def.#1{\tabbox{1}{1}{$#1$}}
\def>#1{\tabbox{2}{1}{$#1$}}
\def^#1{\tabbox{1}{2}{$#1$}}
\def;{\genblankbox{1}{1}\relax}
\catcode`\:=12 \catcode`\.=12 \catcode`\;=12 
\catcode`\>=12 \catcode`\^=7

\newenvironment{tableau}{\bgroup\catcode`\:=13 \catcode`\.=13
  \catcode`\;=13 \catcode`\>=13 \catcode`\^=13
  \setlength{\tabheight}{3ex}\setlength{\tabwidth}{3ex}%
  \def\b##1##2##3{\gentabbox{##1}{##2}{1.2pt}{\vbox{##3}}}%
  \def\n##1##2##3{\gentabbox{##1}{##2}{0.4pt}{\vbox{##3}}}%
  \vbox\bgroup\offinterlineskip}{\egroup\egroup}

\newcommand\Fu{\mathcal{F}_u}

\newcommand\calO{\mathcal{O}}

\newcommand{\FuT}{\mathcal{F}_{u,T}}

\newcommand{\calC}{\mathcal{C}}
\newcommand{\mC}{\mathcal{C}}
\newcommand{\YC}{\mathcal{Y}_\mC}
\newcommand{\YCC}{\mathcal{Y}_{\mC'}}
\newcommand{\OmC}{\overline{\mC}}
\newcommand{\SDTcl}{\Sigma DT_{cl}(\lambda)}
\newcommand{\Irr}{\text{Irr}}

\newtheorem{theorem}{Theorem}[section]
\newtheorem{corollary}[theorem]{Corollary}
\newtheorem{lemma}[theorem]{Lemma}
\newtheorem{proposition}[theorem]{Proposition}

\newtheorem*{theorem*}{Theorem}
\newtheorem*{corollary*}{Corollary}
\newtheorem*{lemma*}{Lemma}

\theoremstyle{definition}
\newtheorem{definition}[theorem]{Definition}
\newtheorem{example}[theorem]{Example}

\newtheorem*{exercise*}{Exercise}
\newtheorem*{question*}{Question}
\newtheorem*{example*}{Example}
\newtheorem*{problem*}{Problem}

\theoremstyle{remark}

\newtheorem*{fact*}{Fact}

\newtheorem*{remark*}{Remark}

\numberwithin{equation}{section}

\begin{document}

\title[Components of the Springer Fiber and Domino Tableaux]{Components of the Springer Fiber and Domino Tableaux}

\author{Thomas Pietraho}
\address{Department of Mathematics, University of Utah, Salt Lake City, UT 84103}
\curraddr{Department of Mathematics, Bowdoin College, Brunswick,
ME 04011}
\email{tpietrah@bowdoin.edu}


\date{}

\begin{abstract}
Consider a complex classical semi-simple Lie group along with the
set of its nilpotent coadjoint orbits.  When the group is of type
$A$, the set of orbital varieties contained in a given nilpotent
orbit is described a set of standard Young tableaux. We
parameterize both, the orbital varieties and the irreducible
components of unipotent varieties in the other classical groups by
sets of standard domino tableaux.  The main tools are
Spaltenstein's results on signed domino tableaux together with
Garfinkle's operations on standard domino tableaux.
\end{abstract}

\maketitle
\section{Introduction}
Let $\mathfrak{g}$ be a  complex semisimple Lie algebra with
adjoint group $G$ and write $\mathcal{O}_f = G \cdot f $ for the
coadjoint orbit of $G$ through $f$ in $\mathfrak{g}^*$.  Fix a
Borel subgroup $B$ of $G$ and let $\mathcal{F}$ be the flag
variety $G/B.$  For a unipotent element $u \in G$,  $\Fu$ is the
variety of flags in $\mathcal{F}$ fixed by the action of $u$. The
orbit $\mathcal{O}_f$ has a natural $G$-invariant symplectic
structure and the Kostant-Kirillov method seeks to attach
representations of $G$ to certain Lagrangian subvarieties of
$\mathcal{O}_f$ (see \cite{gv}, \cite{mihov:thesis}, and
\cite{pietraho:thesis}).  Of particular importance is the set of
orbital varieties, Lagrangian subvarieties of $\calO_f$ that are
fixed by a given Borel subgroup of G.

A result of Spaltenstein identifies the set of orbital varieties
for a given nilpotent orbit with the orbits of a finite group on
the irreducible components of the corresponding unipotent variety
\cite{spaltenstein:fixed}. The main purpose of this paper is to
provide new parameterizations of both, the orbital varieties
contained in a given nilpotent orbit, as well as the irreducible
components of the unipotent variety $\Irr(\Fu)$.

In the case of classical groups, nilpotent coadjoint orbits are
classified by partitions.  Because the number of orbital varieties
contained in a given orbit is finite, one expects that both
orbital varieties and the components of the unipotent variety
should also admit combinatorial descriptions \cite{gerstenhaber}.
This is most apparent when $G$ is of type $A$.

\begin{theorem*}[\cite{spaltenstein:book}] \label{theorem:typea}
Suppose that $G = GL_n(\mathbb{C})$ and the nilpotent orbit
$\mathcal{O}_f$ corresponds to the partition $\lambda$ of $n$.
Then the orbital varieties contained in $\mathcal{O}_f$ as well as
the set of components $\textup{Irr(}\mathcal{F}_u\textup{)}$ are
both parameterized by the family of standard Young tableaux of
shape $\lambda$.
\end{theorem*}

In the setting of other classical groups, a method similar to the
one used to obtain the above can also be employed to describe both
families of objects.  However, the resulting parametrization by
subsets of signed domino tableaux is somewhat cumbersome (see
\cite{spaltenstein:book} and \cite{vanleeuwen:thesis}).  The
following argument suggests a more appealing parameter set.

Let $S$ be the set of partitions indexing the unitary dual of $W$,
the Weyl group of $G$ \cite{mayer} and write $\lambda$ for a
partition lying in $S$.  The dimension of the representation given
by $\lambda$ is precisely the number of standard domino tableaux
of shape $\lambda$.  If we choose a unipotent representative
$u_\lambda \in G$ in the conjugacy class corresponding to
$\lambda$, then Springer's characterization of the representations
$\widehat{W}$ in the top degree cohomology of $\mathcal{F}_u$
\cite{springer78} indicates that
$$
\# SDT(n) = \sum_{\lambda \in S} \dim H^{top}(\mathcal{F}_{u_\lambda},
\mathbb{C}) = \# \{ \Irr(\mathcal{F}_{u_\lambda}) \, | \, \lambda \in S \}
$$

This suggests that $\Irr(\Fu)$ should correspond to a set of
standard domino tableaux in a natural way.  Indeed, this is the
case. The precise relationship between van Leeuwen's parameter set
for $\Irr(\Fu)$ \cite{vanleeuwen:thesis} and the set of domino
tableaux can be described in terms of Garfinkle's notions of
cycles and moving-through maps \cite{garfinkle1}.  After defining
the notion of a {\it distinguished} cycle for a cluster of
dominos, we show that moving through sets of distinguished cycles
of open and closed clusters in van Leeuwen's parameter set defines
a bijection with the set of all domino tableaux of a given size.

\begin{theorem}\label{theorem1}
Suppose that $G$ is a complex classical simple Lie group not of
type $A$.  Then the collection of irreducible components of the
unipotent varieties for $G$ as the unipotent element ranges over
all conjugacy classes is parameterized by $SDT(n)$, the set of
standard domino tableaux of size $n$.
\end{theorem}

The action of the finite group $A_u$ on the irreducible components
$\Irr(\Fu)$ is described by \cite{vanleeuwen:thesis}.  In the
signed domino parametrization, it acts by changing the signs of
open clusters. We exploit this to obtain a parametrization of
orbital varieties by standard domino tableaux. This time, moving
through the distinguished cycles of just the closed clusters in
van Leeuwen's parameter set defines the required bijection.  The
result is a little simpler to state if we consider nilpotent
orbits of the isometry group of a nondegenerate bilinear form,
$G_\epsilon$.

\begin{theorem}\label{theorem2}
Suppose that $G$ is a complex classical simple Lie group not of
type $A$ and $\mathcal{O}$ is the nilpotent orbit of $G_\epsilon$
that corresponds to the partition $\lambda$ . Then the set of
orbital varieties contained in $\mathcal{O}$ is parameterized by
the set of standard domino tableaux of shape $\lambda$.
\end{theorem}

Parameterizations of orbital varieties by domino tableaux have
been obtained in \cite{mcgovern:ssmap}, by describing equivalence
classes in the Weyl group of $G$, as well as in \cite{trapa:rs}.
We will address the compatibility of these parameterizations with
the one above in another paper.

\section{Preliminaries}

We first describe unipotent and orbital varieties, the
relationship between them, and the combinatorial objects we will
use in the rest of the paper.

\subsection{Unipotent and Orbital Varieties}

Let $G$ be a connected complex semisimple algebraic group, $B$ a
Borel subgroup fixed once and for all, and $\mathcal{F} = G/B$ the
flag manifold of $G$. We consider the fixed point set
$\mathcal{F}_u$ of a unipotent transformation $u$ on
$\mathcal{F}$.  It has a natural structure of a projective
algebraic variety, called the {\it unipotent variety}. We write
$\Irr(\mathcal{F}_u)$ for the set of its irreducible components.
The stabilizer $G_u$ of $u$ in $G$ acts on $\mathcal{F}_u$ and
gives an action of its component group $A_u=G_u/G_u^\circ$ on
$\Irr(\mathcal{F}_u).$

Now consider a nilpotent element $f$ of the dual of the Lie
algebra $\mathfrak{g}^*$ of $G$.  Write $\calO^{ad}_f$ for the
orbit of $f$ under the coadjoint action of $G$ on
$\mathfrak{g}^*.$ Using the non-degeneracy of the Killing form, we
can identify $\calO^{ad}_f$ with a subset of $\mathfrak{g}$.  If
$\mathfrak{b}$ is the Lie algebra of $B$ and $\mathfrak{n}$ its
unipotent radical, then the set $\calO^{ad}_f \cap \mathfrak{n}$
inherits the structure of a locally closed algebraic variety from
the orbit $\calO^{ad}_f$.  Its components are Lagrangian
submanifolds of $\calO^{ad}_f$ and are known as {\it orbital
varieties} \cite{ginsburg}.  There is a simple relationship
between the set of orbital varieties contained in a given
nilpotent orbit and the irreducible components of the
corresponding unipotent variety. Suppose that the unipotent
element $u$ of $G$ and the nilpotent element $f$ of
$\mathfrak{g}^*$ correspond to the same partition.

\begin{theorem}[\cite{spaltenstein:fixed}]
There is a natural bijection
$$ \textup{Irr}(\calO^{ad}_f \cap \mathfrak{n}) \longrightarrow \textup{Irr}(\mathcal{F}_u)/A_u$$
between the orbital varieties contained in the nilpotent orbit
$\calO^{ad}_f$ and the orbits of the finite group $A_u$ on
$Irr(\mathcal{F}_u).$
\end{theorem}

The set of nilpotent orbits for a classical $G$ admits a
combinatorial description by partitions. Write $\mathcal{P}(n)$
for the set of partitions $\lambda = [\lambda_1, \lambda_2,
\ldots, \lambda_k]$ of $n$, ordered so that $\lambda_i \geq
\lambda_{i+1}$.

\begin{theorem} Nilpotent orbits in $\mathfrak{gl}_n$ are in one-to-one
correspondence with the set $\mathcal{P}(n).$
\end{theorem}

The corresponding statement for the other classical groups is not
much more difficult.  To obtain slightly cleaner statements, we
will state it in terms of the nilpotent orbits of the sightly
larger isometry groups of  nondegenerate bilinear forms. Let
$\epsilon = \pm 1$ and consider a nondegenerate bilinear form on
$\mathbb{C}^m$ satisfying $(x,y)_\epsilon = \epsilon
(y,x)_\epsilon$ for all $x$ and $y$. Let $G_\epsilon$ be the
isometry group of this form and $\mathfrak{g}_\epsilon$ be its Lie
algebra. Define a subset $\mathcal{P}_{\epsilon}(m)$ of
$\mathcal{P}(m)$ as the partitions $\lambda$ satisfying $
\#\{j|\lambda_j=i\}$ is even for all $i$ with $(-1)^i=\epsilon$.
The classification of nilpotent orbits now takes the form:

\begin{theorem}[\cite{gerstenhaber}] Let $m$ be the dimension of the
standard representation of  $G_\epsilon$. Nilpotent
$G_\epsilon$-orbits in $\mathfrak{g}_\epsilon$ are in one to one
correspondence with the partitions of $m$ contained in
$\mathcal{P}_{\epsilon}(m).$
\end{theorem}

The nilpotent $G_\epsilon$ orbits in $\mathfrak{g}_\epsilon$ can
be identified with the nilpotent orbits of the corresponding
adjoint group with one exception. In type $D$, precisely two
nilpotent orbits of the adjoint group correspond to every very
even partition. We will write $\mathcal{O}_f$ for the
$G_\epsilon$-orbit through the nilpotent element $f$ and
$\mathcal{O}_\lambda$ for the $G_\epsilon$-orbit that corresponds
to the partition $\lambda$ in this manner.

The group $A_u$ is always finite, and in the setting of classical
groups, it is always a two-group.  More precisely:

\begin{theorem}[\cite{vanleeuwen:thesis}(2.4.1)]
The group $A_u$ is always trivial when $G$ is of type $A$.  In the
other classical types, let $B_{\lambda}$ be the set of the
distinct parts of $\lambda$ satisfying  $(-1)^{\lambda_i} = -
\epsilon$. Then $A_u$ is a \textup{2}-group with $|B_{\lambda}|$
components.
\end{theorem}

\subsection{Standard Tableaux}
A partition of of an integer $m$ corresponds naturally to a Young
diagram consisting of $m$ squares.  We call the partition
underlying a Young diagram its {\it shape}.  Recall the
definitions of the sets of standard Young tableaux and standard
domino tableaux from, for instance, \cite{garfinkle1}. We will
write $SYT(\lambda)$ and $SDT(\lambda)$ respectively for the sets
of Young and domino tableaux of shape $\lambda$. We refer to both
objects generically as {\it standard tableaux} of shape $\lambda$,
or $ST(\lambda),$ hoping that the precise meaning will be clear
from the context. Also, we will write $ST(n)$ for the set of all
standard tableaux with largest label $n$.

We view each standard tableau $T$ as a set of ordered pairs
$(k,S_{ij})$, denoting that the square in row $i$ and column $j$
of $T$ is labelled by the integer $k$. When $T$ is a domino
tableau, the domino with label $k$, or $D(k,T)$, is a subset of
$T$ of the form  $\{(k,S_{ij}),(k,S_{i+1,j})\}$ or
$\{(k,S_{ij}),(k,S_{i,j+1})\}$.  We call these vertical and
horizontal dominos, respectively.  For convenience, we will refer
to the set $\{(0,S_{11})\}$  as the zero domino when in type $B$.
Whenever possible, we will omit labels of the squares and write
$S_{ij}$ for $(k,S_{ij})$.  In that case, define $label \, S_{ij}
= k$.

\begin{definition}
For a standard tableau $T$, let $T(i)$ denote the tableau formed
by the squares of $T$ with labels less than or equal to $i$. A
domino tableau $T$ is {\it admissible} of type $X$ = $B$, $C$, or
$D$, if  the shape of each $T(i)$ is a partition of a nilpotent
orbit of type $X$.
\end{definition}

The dominos that appear within admissible tableaux fall into three
categories. Following \cite{vanleeuwen:thesis} , we call these
types  $I^+,$ $I^-,$ and $N$.
\begin{example}
Suppose that $G$ is of type $C$ and consider the tableaux

$$
\raisebox{3ex}{$T=$ \;}
\begin{tableau}
:^1^2^3^4>5\\
:;\\
\end{tableau}
\hspace{1in}
\raisebox{3ex}{$T'=$ \;}
\begin{tableau}
:^1>2>3>5\\
:;>4\\
\end{tableau}
$$

\noindent Then $T$ is admissible of type $C$ but $T'$ is not,
since $shape \, T'(2) = [3,1]$ is not the partition of a nilpotent
orbit in type $C$.  The dominos $D(1,T)$ and $D(3,T)$ are of type
$I^-$, $D(2,T)$ and $D(4,T)$ are of type $I^+,$ and $D(5,T)$ is of
type $N$.
\end{example}

We also recall the notions of a cycle in a domino  tableau and
moving through such a cycle, as defined in \cite{garfinkle1}. We
will think of cycles as both, subsets of dominos of $T$, as well
as just sets of their labels. Write $MT(D(k,T),T)$ for the image
of the domino $D(k,T)$ under the moving through map and $MT(k,T)$
for the image of $T$ under moving through the cycle containing the
label $k$. If $U$ is a set of cycles of $T$ that can be moved
through independent of one another, we will further abuse notation
by writing $MT(U,T)$ for the tableau obtained by moving through
all the cycles in $U$.  Recall the definition of $X$-fixed and
$X$-variable squares for $X$ = $B$, $C$, $D$, or $D'$
\cite{garfinkle1}. Under the moving through map, the labels of the
fixed squares are preserved while those of variable ones may
change.  We will call a cycle whose fixed squares are $X$-fixed an
$X$-{\it cycle}.  Note also that the $B$- and $C$-cycles as well
as the $D$- and $D'$-cycles in a given tableau $T$ coincide.

\begin{example}
Consider the domino tableaux $T$ and $T'$ from the previous
example. The $C$-cycles in $T$ are \{1\}, \{2,3\}, and \{4,5\}
while those in $T'$ are \{1\} and \{2,3,4,5\}.  We have

$$
\raisebox{3ex}{$MT(2,T)=$ \;}
\begin{tableau}
:^1>2^4>5\\
:;>3\\
\end{tableau}
\hspace{.3in}
\raisebox{3ex}{$MT(4,T)=$ \;}
\begin{tableau}
:^1^2^3>4>5\\
:;\\
\end{tableau}
$$

\noindent The $D$-cycles in $T$ are \{1,2\}, \{3,4\}, and \{5\},
while there is only one in $T'$, mainly \{1,2,3,4,5\}.
\end{example}

\section{Signed Domino Tableaux Parameterizations}

The irreducible components of the unipotent variety $\Fu$ for
classical $G$ were described by N. Spaltenstein in
\cite{spaltenstein:book}.  We summarize this parametrization as
interpreted by M.A. van Leeuwen \cite{vanleeuwen:thesis}.  Its
advantage lies in a particularly translucent realization of the
action of $A_u$ on $\Irr(\mathcal{F}_u)$.

Fix a unipotent element $u \in G$ and let $\lambda_u$ be the
partition of the corresponding nilpotent orbit.  Define a map
$$
\Gamma_u: \mathcal{F}_u \longrightarrow ST(\lambda_u)
$$
by the following procedure.  Fix a flag $F \in \mathcal{F}_u$.
Adopting  notation of \cite{vanleeuwen:thesis}, let
$\lambda^{(i)}$ be the shape of the Jordan form of the unipotent
operator induced by $u$ on $F^{(i)}$. The difference between the
Young diagrams $\lambda^{(i)}$ and $\lambda^{(i+1)}$ is  one
square in type A and a domino in the other classical types
\cite{spaltenstein:book}.  By assigning the label $i+1$ to the set
$\lambda^{(i+1)} \setminus \lambda^{(i)}$ for each $i$, we obtain
 a standard tableau of shape $\lambda_u.$

\begin{theorem} When $G$ is of type $A$, the map $\Gamma_u$ defines a
surjection onto $SYT(\lambda_u)$ that separates points of
$\textup{Irr}(\mathcal{F}_u)$.  That is, it defines a bijection
$$
\Gamma_u: \textup{Irr}(\mathcal{F}_u) \longrightarrow SYT(\lambda_u).
$$
\end{theorem}

\begin{corollary} When $G$ is of type $A$, the orbital
varieties $\textup{Irr}(\mathcal{O}_\lambda \cap \mathfrak{n})$
are parameterized by  the set $SYT(\lambda)$.
\end{corollary}

In the other classical types, any domino tableau in the image of
$\Gamma_u$ is admissible.  Admissible tableaux, however, do not
fully separate the components of $\Fu$.  If two flags give rise to
different domino tableaux in this way, they lie in different
components of $\Fu$.  However, the converse is not true.  The
inverse image $\FuT$ of a given admissible tableau $T$ under this
identification is in general not connected. Nevertheless, the
irreducible components of $\FuT$ are precisely its connected
components \cite{vanleeuwen:thesis}(3.2.3).  Accounting for this
disconnectedness yields a parametrization of $\Irr(\Fu)$.

\begin{definition}
A signed domino tableau $T$ of shape $\lambda$ is an admissible
domino of shape $\lambda$ with a choice of sign for each domino
of type $I^+$.  The set of signed domino tableaux is denoted
$\Sigma DT(\lambda)$.
\end{definition}

The set $\Sigma DT(\lambda_u)$ is too large to parameterize
$\Irr(\mathcal{F}_u)$ and we follow \cite{vanleeuwen:thesis} in
defining equivalence classes.  We recall the notion of a {\it
cluster} of dominos.

\begin{definition} \label{definition:openclosed}
  Write $cl(0)$ for the cluster containing $D(1,T)$ in types $B$ and $C$. A cluster is {\it open} if it contains an $I^+$ or $N$ domino along its right edge and is not $cl(0)$.  A cluster that is neither $cl(0)$ nor open is {\it closed}.  Denote the set of open clusters of $T$ by $OC(T)$ and the set of closed clusters as $CC(T)$.  For a cluster $\mathcal{C}$, let $I_\mathcal{C}$ be the domino in $\mathcal{C}$ with the smallest label and take $S_{ij}$ as its left and uppermost square.  For $X$ equal to $B$ or $C$, we say that $\mathcal{C}$ is an $X$-cluster iff $i+j$ is odd.  For $X$ equal to $D$ or $D'$, we say that $\mathcal{C}$ is an $X$-cluster iff $i+j$ is even.
\end{definition}

The definition of open and closed differs from \cite{vanleeuwen:thesis} as we do not
call $cl(0)$ an open cluster.   The open clusters of $T$
correspond to the parts of $\lambda$ contained in $B_{\lambda}$,
the set parameterizing the $\mathbb{Z}_2$ factors of
$A_{\lambda}.$

\begin{definition}
If $T, T' \in \Sigma DT(\lambda)$, let $T \sim_{op,cl} T'$ iff
the underlying tableaux are the same and the products of signs in
all corresponding open and closed clusters of $T$ and $T'$ agree.
Denote the equivalence classes by $\Sigma DT_{op,cl}(\lambda).$
Define the set $\Sigma DT_{cl}(\lambda)$ similarly.  We represent
the elements of $\Sigma DT_{op,cl}(\lambda)$ and $\Sigma
DT_{cl}(\lambda)$ as admissible tableaux with a choice of sign for
each of the appropriate clusters.
\end{definition}

There is a considerable amount of freedom in how a flag of
$\mathcal{F}_{u,T}$ can be assigned an equivalence class of signed
admissible domino tableaux.  A particular choice is presented in
\cite{vanleeuwen:thesis}(3.4), defining a map
$$\widetilde{\Gamma}_u: \mathcal{F}_u \longrightarrow \Sigma
DT_{op,cl} (\lambda_u).$$

We describe an action of $A_{u}$ on $\Sigma
DT_{op,cl}(\lambda_u).$ For $r \in B_\lambda$, let $b_T(r)$ be the
cluster that contains a domino ending a row of length $r$ in $T$.
Let $\xi_{r}$  act trivially if $b_{T}(r) = cl(0)$ and by changing
the sign of the open cluster $b_T(r)$ otherwise. For each $r \in
B_{\lambda}$, let $g_r$ denote the generator of the corresponding
$\mathbb{Z}_2$ factor of $A_u$. One can now define the action of
$g_r$ on $\Sigma DT_{op,cl}(\lambda_u)$ by $g_r[T] = \xi_r[T].$

\begin{theorem}[\cite{vanleeuwen:thesis}]
Suppose that $G$ is a classical group not of type A and $u$ is a
unipotent element of $G$ corresponding to the partition $\lambda$
. The map $\widetilde{\Gamma}_u$ defines an $A_{u}$-equivariant
bijection between the components $\textup{Irr}(\Fu)$ and $\Sigma
DT_{op,cl}(\lambda).$
\end{theorem}

Since $A_u$ acts by changing the signs of the open clusters of
$\Sigma DT_{op,cl}(\lambda),$ it is simple to parameterize the
$A_u$ orbits on $\Irr(\Fu).$

\begin{corollary}\label{theorem:oun}
Suppose that $G$ is a classical group not of type A and
$\mathcal{O}'_\lambda$ is the nilpotent orbit corresponding to the
partition $\lambda$ . The orbital varieties
$\textup{Irr}(\mathcal{O}_\lambda \cap \mathfrak{n})$ are
parameterized by $\Sigma DT_{cl}(\lambda).$
\end{corollary}

\section{Domino Tableaux Parameterizations}

We show how to index the components $\Irr(\Fu)$ and
$\Irr(\mathcal{O}_\lambda \cap \mathfrak{n})$ by families of
standard tableaux. In type $A$, this is Theorem
\ref{theorem:typea}. For the other classical types, we define maps
from domino tableaux with signed clusters to the set of standard
domino tableaux by applying Garfinkle's moving through map to
certain distinguished cycles.

\subsection{Definition of Bijections}

Consider an $X$-cluster $\mC$ and let $I_C$ be the domino in $\mC$
with the smallest label.  Let $\YC$ be the X-cycle through $I_C$.
We call it  the {\it initial} cycle of the cluster $\mC$.

\begin{proposition}
\label{proposition:cycleincluster} A cluster of an admissible
domino tableau $T$ that is either  open or closed contains its
initial cycle.
\end{proposition}

We defer the proof to another section.  Armed with this fact, we can propose a map
$$\Phi: \Sigma DT_{op,cl}(n) \longrightarrow SDT(n)$$
by moving through the distinguished cycles of all open and closed
clusters with positive sign. More explicitly, for a tableau $T \in
\Sigma DT_{op,cl},$ let $C^+(T)$ denote the set of open and closed
clusters of $T$ labelled by a $(+)$ and let $\sigma (T) = \{ \YC
\; | \; \mC \in C^+(T) \}$ be the set of their distinguished
cycles.  Write $|T|$ for the standard domino tableau underlying $T$.  We define $$\Phi(T) = MT (\sigma(T), |T|).$$

\begin{lemma}
The map $\Phi : \Sigma DT_{op,cl}(n) \longrightarrow SDT(n)$  is a
bijection. We can view the set $\Sigma DT_{cl}(n)$ as a subset of
$ \Sigma DT_{op,cl}(n)$ by assigning a negative sign to each
unsigned open cluster of a domino tableau in $\Sigma DT_{cl}(n).$
Restricted to $\Sigma DT_{cl}(n)$, $\Phi$  preserves the shapes of
tableaux and defines a bijection $\Phi: \Sigma DT_{cl}(\lambda)
\longrightarrow SDT(\lambda)$ for each $\lambda$ a shape of a nilpotent orbit.
\end{lemma}

\begin{proof}
We check that $\Phi$ is well-defined, that its image lies in
$SDT(n)$, and then construct its inverse.  We first need to know that
the definition of $\Phi$ does not depend on which order we move
through the cycles in $\sigma(T)$. It is enough to check that if
$\YC$ and $\YCC \in \sigma(T)$, then $\YCC$ is also lies in
$\sigma(MT(|T|, \YC))$. While this statement is not true for
arbitrary cycles, in our setting, this is Lemma \ref{lemma:well}.

The image of $\Phi$ indeed lies in $SDT(n)$. That
$\Phi(T)$ is itself a domino tableau follows from the fact that
moving through any cycle of $|T|$ yields a domino tableau.   Hence $\Phi(T) \in SDT(n)$ and if $T \in \Sigma DT_{cl}(\lambda)$ then $\Phi(T) \in SDT(\lambda)$ since in this case $\Phi$ moves through only closed cycles.

The definition of a cluster forces the initial
domino $I_{\mC}$ of every closed cluster to be of type $I^+$.
By the definition of moving through, the image of
$MT(I_{\mC},T)$ in $MT(\YC,T)$ is inadmissible, i.e. it is a
horizontal domino not of type $N$.  In general, all the
inadmissible dominos in $\Phi(T)$ appear within the image of
distinguished cycles under moving through.  Furthermore, the
lowest- numbered domino within each cycle is the image of the
initial domino of some distinguished cycle.  With  this
observation, we can construct the inverse of $\Phi.$ We define a
map $$\Psi : \Phi (\Sigma DT_{op,cl}(n)) \longrightarrow \Sigma DT_{op,cl}(n)$$ that satisfies
$\Psi \circ \Phi = \text{Identity}$.  Let $\iota(\Phi (T))$ be the
set of cycles in $\Phi(T)$ that contain inadmissible dominos.  We
define $\Psi (\Phi (T)) = MT (\Phi (T), \iota (\Phi (T))).$ By the
above discussion, $\iota(\Phi (T))$ contains precisely the images
of cycles in $\sigma (T)$.  Hence

$$\Psi (\Phi (T)) = MT(\Phi(T),\iota(\Phi(T)))=MT(MT(|T|,\sigma(T)))=T$$
as desired.  Thus $\Phi$ is a bijection onto its image in
$SDT(n)$ and restricted to $\Sigma DT_{cl}(\lambda)$, it is a bijection with its image in $SDT(\lambda)$.  As we already know that the sets $\SDTcl$ and
$SDT(\lambda)$ both parameterize the same set of orbital varieties,
and that $\Sigma DT_{op,cl}(n)$ and $SDT(n)$ both parameterize the same set of irreducible components of unipotent varieties,
$\Phi$ must provide bijections between these two sets.
\end{proof}

Theorems \ref{theorem1} and \ref{theorem2} are immediate consequences.

\begin{example}
Let $G$ be of type $D$ and suppose that both $u$ and
$\mathcal{O}_\lambda$ correspond to the partition $\lambda =
[3^2]$. The van Leeuwen parameter set $\Sigma DT_{op,cl}([3^2])$
for $\Irr(\mathcal{F}_u)$ is:
$$
\begin{small}
\begin{tableau}
:^{1 \atop +}^2^{3 \atop +}\\
:;\\
\end{tableau}
\hspace{.2in}
\begin{tableau}
:^{1 \atop -}^2^{3 \atop +}\\
:;\\
\end{tableau}
\hspace{.2in}
\begin{tableau}
:^{1 \atop +}^2^{3 \atop -}\\
:;\\
\end{tableau}
\hspace{.2in}
\begin{tableau}
:^{1 \atop -}^2^{3 \atop -}\\
:;\\
\end{tableau}
\hspace{.2in}
\begin{tableau}
:^{1 \atop +}>2\\
:;>3\\
\end{tableau}
\hspace{.2in}
\begin{tableau}
:^{1 \atop -}>2\\
:;>3\\
\end{tableau}
\end{small}
$$
The image of  $\Sigma DT_{op,cl}([3^2])$ under $\Phi$ is the
following set of standard domino tableaux.  We write the image of
a given tableau in the same relative position.  Note that this
parameter set for $\Irr(\mathcal{F}_u)$ consists of all tableaux
of shapes $[3^2]$ and $[4,2]$.
$$
\begin{small}
\begin{tableau}
:>1>3\\
:>2\\
\end{tableau}
\hspace{.1in}
\begin{tableau}
:^1^2>3\\
:;\\
\end{tableau}
\hspace{.1in}
\begin{tableau}
:>1^3\\
:>2\\
\end{tableau}
\hspace{.1in}
\begin{tableau}
:^1^2^3\\
:;\\
\end{tableau}
\hspace{.1in}
\begin{tableau}
:>1>2\\
:>3\\
\end{tableau}
\hspace{.1in}
\begin{tableau}
:^1>2\\
:;>3\\
\end{tableau}
\end{small}
$$
The van Leeuwen parameter set $\Sigma DT_{cl}([3^2])$ for the
orbital varieties contained in $\mathcal{O}_\lambda$ is:
$$
\begin{small}
\begin{tableau}
:^{1 \atop +}^2^3\\
:;\\
\end{tableau}
\hspace{.8in}
\begin{tableau}
:^{1 \atop -}^2^3\\
:;\\
\end{tableau}
\hspace{.8in}
\begin{tableau}
:^1>2\\
:;>3\\
\end{tableau}
\end{small}
$$
Its image under $\Phi$ is the set of all domino tableaux of shape
$[3^2]$.  Again, we write the image of a tableau in the same
relative position.
$$
\begin{small}
\begin{tableau}
:>1^3\\
:>2\\
\end{tableau}
\hspace{.8in}
\begin{tableau}
:^1^2^3\\
:;\\
\end{tableau}
\hspace{.8in}
\begin{tableau}
:^1>2\\
:;>3\\
\end{tableau}
\end{small}
$$
\end{example}
\subsection{Independence of Moving Through Initial Cycles}

\begin{lemma} Consider open or closed clusters $\calC$ and $\calC'$ and
 their initial cycles $\YC$ and $\YCC.$  Then $\YC$ is again a cycle
 in $MT(|T|, \YCC)$. \label{lemma:well}
\end{lemma}

\begin{proof}
If $\calC$ and $\calC'$ are clusters of the same type, then so are
their initial cycles and the lemma is \cite{garfinkle1}(1.5.29).
Otherwise, without loss of generality, take $\calC$ to be a
$C$-cluster and $\calC'$ to be a $D$-cluster.  As the proof in the
other cases is similar, we can also assume that $\YC$ is $C$-boxed
while $\YCC$ is $D$-boxed.

Suppose that the dominos $D(r) \in \YC$ and $D(s) \in \YCC$ lie in
relative positions compatible with the diagram

$$
\begin{tiny}
\begin{tableau}
:.{}.s\\
:.r\\
\end{tableau}
\end{tiny}
$$
where the box labelled by $r$ is fixed.  The same squares in
$MT(|T|,\YCC)$ have the labels

$$
\begin{tiny}
\begin{tableau}
:.{}.{s'}\\
:.r\\
\end{tableau}
\end{tiny}
$$
for some $s'$.

To prove the lemma, we need to show that $s<r$ implies $s'<r$ and
$s>r$ implies $s'>r.$ Since our choice of $r$ and $s$ was
arbitrary, this will show that $\YC$ remains a cycle. There are
two possibilities for the domino $D(s)$.  It is either horizontal
or vertical and must occupy the following squares:

\begin{center}
\begin{tiny}
$
\begin{tableau}
:.s.s\\
:.r\\
\end{tableau}
$ \hspace{1in}
$
\begin{tableau}
:.{}.s\\
:.{}.s\\
:.r\\
\end{tableau}
$
\end{tiny}
\end{center}

\begin{center}

Case (i) \hspace{.8in} Case (ii)

\end{center}

Case (i).  In this case, $s<r$ always.  Garfinkle's rules for
moving through  imply that $MT(|T|, D(r)) \cap \calC' \neq
\varnothing.$  This is a contradiction since we know by hypothesis
that $\YC \neq \YCC$ .  Hence this case does not occur.

Case (ii).  First suppose $s>r$.  Then the our squares within $MT(|T|, \YCC)$ must look like

$$
\begin{tiny}
\begin{tableau}
:.{}.s\\
:.{}.{s'}\\
:.r\\
\end{tableau}
\end{tiny}
$$
for some $s' \neq s.$  Since  the tableau $MT(\YCC,T)$ is
standard, this requires that  $s'>s$ implying $s'>r$ which is what
we desired.  Now suppose $s<r$ and suppose the squares in our
diagram look like
$$
\begin{tiny}
\begin{tableau}
:.{}.s\\
:.t.s\\
:.r.u\\
\end{tableau}
\end{tiny}
$$
As in Case (i), we find that  $D(t) \notin \calC'$.  Since $D(t)
\in \calC$, type $D(s) = I^+$ implies type $D(t) = I^-$, type
$D(r) = I^-,$ and type $D(u) = I^+$.  Otherwise, the rules
defining clusters would force $s$ to lie in the cluster  $\calC$.
Now $D(u)$ lies in the initial cycle of a closed cluster of same
type as $\calC'$.  Since it lies on the periphery and its type is
$I^+$, then its top square must be fixed.  In particular, $D(u)
\notin \calC$. But $s<r$ implies $MT(D(r)) \cap D(u) \neq
\varnothing$.  This is a contradiction, implying that this case
does not arise.

To finish the proof, we must examine the possibility that $D(s)$
and $D(r)$ lie in the relative positions described by
$$
\begin{tiny}
\begin{tableau}
:.{}.r\\
:.s\\
\end{tableau}
.
\end{tiny}
$$
This case is completely analogous and we omit the proof.
\end{proof}

This lemma shows that the image of moving though a subset of
distinguished cycles is independent of the order in which these
cycles are moved though.  Note, however, that a similar result is
not true for subsets of arbitrary cycles.
\subsection{Nested Clusters and the Periphery of a Cluster}

We aim to show that closed and open clusters contain their
distinguished cycles.  The proof has two parts.  First, we show
that $\YC$ is contained in a larger set of clusters $\OmC$,
defined as the union of $\mathcal{C}$ with all of its {\it nested}
clusters. Then, we show that $\YC$ intersects each of the nested
clusters trivially.

Let $\mathcal{C}$
be a cluster of a tableau $T$ and denote by $row_k
\,T=\{S_{k,j}|\; j \geq 0\}$ the $k$th row of $T$.  Define
$col_k\,T$ similarly.  If $row_k\,T \cap \mathcal{C} \neq
\emptyset$, let $\inf_k\,\mC = \inf\{j|\; S_{k,j} \in row_k\,T
\cap \mC \}$ and $\sup_k\,\mC = \sup\{j\;|\; S_{k,j} \in row_kT
\cap \mC\} $.

\begin{example}
Consider the following tableau of type $D.$  It has two closed
clusters given by  the sets
$\mathcal{C}=\{1,2,3,4,5,8,9,10,11,12\}$ and $\mathcal{C'}=
\{6,7\}$.
$$
\begin{tiny}
\begin{tableau}
:^1>3>5^{11}\\
:;^4^6^7^8\\
:^2;;;;^{12}\\
:;>9>{10};\\
\end{tableau}
\end{tiny}
$$
$\mathcal{C}$ is a $D$-cluster while $\mathcal{C'}$ is a
$B$-cluster.  $\YC$ is then a $D$-cycle
and consists of the dominos in the set $\{1,3,5,11,12,10,9,2\}.$  $T$ has two other
$D$-cycles, $\{4,6\}$ and $\{7,8\}.$ Both intersect
$\mathcal{C}$, but are not contained within it.  The $B$-cycle
$\YCC$ equals $ \{6,7\} $ and is contained
in $\mathcal{C'}$. Hence an $X$-cluster may not contain all the
X-cycles through its dominos.  However, it always contains its
initial cycle.  Also notice that $\mathcal{C}$ completely
surrounds $\mathcal{C'}$.  We call such interior clusters {\it
 nested}.
\end{example}

Nested clusters complicate the description of clusters.  To
simplify our initial results, we would like to consider the set formed
by a cluster together with all of its nested clusters.  To be more precise:

\begin{definition}
Let $\mC'$ be a cluster of $T$.  It is {\it nested} in $\mC$ if
all of the following are satisfied:
\begin{align*}
        \inf\{k|row_kT \cap \mC' \neq \emptyset\} &> \inf\{k|row_kT \cap \mC \neq
        \emptyset\}\\
        \sup\{k|row_kT \cap \mC' \neq \emptyset\} &< \sup\{k|row_kT \cap \mC \neq
        \emptyset\}\\
        \inf\{k|col_kT \cap \mC' \neq \emptyset\} &> \inf\{k|col_kT \cap \mC \neq
        \emptyset\}\\
        \sup\{k|col_kT \cap \mC' \neq \emptyset\} &< \sup\{k|col_kT \cap \mC \neq
        \emptyset\}\\
\end{align*}
Define $\overline {\mC}$ to be the union of $\mC$ together with
all clusters nested within it.  We will write
$periphery(\overline{\mC})$ for the set of dominos in $\overline
{\mC}$ that are adjacent to some square of $T$ that does not lie
in $\overline {\mC}$. Note that $periphery(\OmC)$ is a subset of
the original cluster $\mC$.
\end{definition}

\begin{example}
In the above tableau, $\mathcal{C'}$ is nested in
$\mathcal{C}$.  Furthermore, $\mathcal{C} \cup \mathcal{C'} =
\OmC =T,$ and $periphery (\OmC) = \mathcal{Y}_{\mathcal{C}}
\subset $ $\mathcal{C}$.
\end{example}

The next two propositions describe properties of dominos that
occur along the left and right edges of $\OmC$.   Recall that our
definition of the cycle $\YC$ endows $\mC$ as well as $\OmC$ with
a choice of fixed and variable squares by defining the left and
uppermost square of $I_{\mC}$ as fixed.

\begin{proposition} \label{proposition:infsup}
Suppose that $\mC$ is a non-zero cluster of a domino tableau $T$
and that the intersection of the $k$-th row of $T$ with $\mC$ is
not empty. Then the dominos $D(label(T_{k,inf_k\,\mC}),T)$ and
$D(label(T_{k,inf_k\,\OmC}),T)$ are both of type $I^+.$ If $\mC$
is also closed, then the dominos $D(label(T_{k,sup_k\,\mC}),T)$
and $D(label(T_{k,sup_k\,\OmC}),T)$ are of type $I^-.$
\end{proposition}
\begin{proof}
The first statement is true for all non-zero clusters by
\cite{vanleeuwen:thesis}(3.3).  The second statement is the
defining property of closed clusters.
\end{proof}

\begin{proposition}
\label{theorem:periphery} Suppose that $\mC$ is a non-zero cluster
of a domino tableau $T$.  If the domino $D$ consisting of the
squares  $S_{pq}$ and  $S_{p+1,q}$  lies in $periphery(\OmC) $,
then
\begin{enumerate}
\item $ S_{pq}$ is fixed if $type$ $D = I^+$ and
\item $S_{p+1,q}$ is fixed if $type$ $D = I^-$
\end{enumerate}
\end{proposition}
\begin{proof}
Case $(i)$.  Assume that there is a $D'$ in the $periphery(\OmC)$ of
type $I^+$ whose uppermost square is not fixed.  Then
$periphery(\OmC)$ must contain two type  $I^+$ dominos $E =
\{S_{kl}, S_{k+1,l} \}$ and $E' = \{S_{k+1,m}, S_{k+2,m} \}$ with
the squares $S_{kl}$ and  $S_{k+2,m}$ fixed and $|m-l|$ minimal.

Assume $m<l$. The opposite case can be proved by a similar
argument.  Because $E'$ is of type $I^+$, there is an integer $t$
such that $m<t<l$, $S_{k+1,t} \in periphery(\OmC)$, and $t$ is
maximal with these properties.  Let $F$ be the domino containing
$S_{k+1,t}$.  $F$ has to be  $\{S_{k+1,t}, S_{k+2,t} \}$ and of
type $I^-$.  If its type was $I^-$ or $N$,
\cite{vanleeuwen:thesis}(3.3 (17)) would force $S_{k+1,t+1}$ to be
in $periphery(\OmC)$ as well.  If $F$ on the other hand was
$\{S_{k+1,t}, S_{k,t} \}$, this would contradict the minimality of
$|m-l|$.  We now consider two cases.

\begin{itemize}
\item[(a)] Assume $t=l-1$.  Because $E$ and $F$ lie in
$periphery(\OmC)$ and hence in $\mC$, $\mC$ must contain a domino
of type $N$ of the form $\{S_{u,l-1}, S_{u,l} \}$ with $u>k+2$
and $u$ minimal with this property.  The set of squares
$\{S_{p,l-1} | k+2 <p<u \} \cup \{S_{pl} | k+1<p<u \}$ must be
tiled by dominos, which is impossible, as its cardinality is odd.

\item[(b)] Assume $t<l-1$.  We will contradict the maximality of
$t$.  Because $E$ and $F$ both lie in $\mC$, $\mC$ must contain a
sequence $H_{\alpha}$ of dominos of type $N$ satisfying
$$H_\alpha =
\{S_{k+1+f(\alpha),t+2\alpha},S_{k+1+f(\alpha),t+2\alpha +1} \}$$
where $0\leq \alpha \leq \frac{l-t+1}{2}$.  We choose each
$H_\alpha$ such that for all $\alpha,$ $f(\alpha)$ is minimal with
this property.  Because the sets $ \{ S_{k+p,l} | k+1
<p<k+1+f(\frac{l-t+1}{2}) \}$ and  $ \{ S_{k+p,t} | k+2
<p<k+1+f(0) \}$ have to be tiled by dominos of type $I^+$ and
$I^-$ respectively, $f(0)$ has to be even and
$f(\frac{l-t+1}{2})$ has to be odd.  Hence there is a $\beta$ such
that $f(\beta)$ is even and $f(\beta +1)$ is odd.

Assume  $f(\beta) < f(\beta +1)$, but the argument in the other
case is symmetric.  Let $G$ be the domino containing the square $
S_{k+1+f(\beta),t+2\beta +2}.$  $G$ must belong to $\mC$, as
$H_{\beta}$ and $G$ is either of type $I^-$ or $N$.  The type
of $G$ cannot be $N$, however, as this would contradict the
condition on $f$.  Hence $G$ must be of type $I^-$.  If $G$
equals  $\{ S_{k+1+f(\beta),t+2\beta+2},S_{k+f(\beta),t+2\beta +2}
\}$.  Then by successive applications of
\cite{vanleeuwen:thesis}(3.3 (17)), the set of dominos
        $$\{ \{S_{k+f(\beta) - \gamma \epsilon,t+2\beta + \epsilon},S_{k+1+f(\beta) - \gamma - \epsilon , t + 2\beta + \epsilon}\}\}$$ with  $\epsilon = 1 \text{ or }2 \text{ and }0 \leq \gamma \leq f(\beta) -2 $
is contained in $\mC$ as well.  But this means that
$t+2\beta+\epsilon$ for $\epsilon = 1 \text{ or } 2$ satisfies the
defining property of $t$, contradicting its maximality.
\end{itemize}
Case $(ii)$. We would like to show that the bottom square is fixed for
every $I^-$ domino in $periphery (\OmC)$.  It is enough to show
that this is true for one such domino, as an argument similar to
that in case $(i)$  can be repeated for the others.  Let $l= \inf
\{k|row_kT \cap \OmC = \emptyset \}$.  Then by
\ref{proposition:infsup} and the definition of fixed, we know that
$S_{l,\inf_l \OmC}$ is fixed.  As
$\{S_{l,\sup_l\OmC},S_{l+1,\sup_l \OmC} \}$ is a domino of type
$I^-$ in $periphery(\OmC)$, we have found the desired domino.
\end{proof}

\begin{lemma}
\label{lemma:cycleinclusterii} The following inclusions hold when
$\mathcal{C}$ is an open or closed cluster: $periphery \, (\OmC)
\subset \YC \subset \overline{\mC}$.
\end{lemma}
\begin{proof}
Recall that our choice of a fixed square in $I_{\mathcal{C}}$
defines the fixed squares in all of $\OmC$.  Define
$\widetilde{\mathcal{C}}$ as $\OmC$ when $\mC$ is closed and
$\OmC$ union with all empty holes and corners of $|T|$ adjacent to
$\mC$ when $\mC$ is open \cite{garfinkle1}(1.5.5).  We show that
the image $MT(D,T)$ of $D$ in $periphery(\OmC)$ lies in
$\widetilde{\mC}$.  This shows the second inclusion, as if any
domino in $periphery(\OmC)$ stays in $\OmC$ under moving through,
then so must the cycle $\YC$.  The first inclusion is a
consequence of the argument and the definitions of moving through
and clusters. We differentiate cases accounting for different
domino positions along $periphery(\OmC).$

Case $(i)$.  Take $D=\{(k,S_{ij}),(k,S_{i+1,j})\}$ and suppose $type \, D
  = I^+$.  Because $D$ lies on $periphery (\OmC)$,
  Proposition \ref{theorem:periphery} implies that $S_{ij}$ is
  fixed.  Due to \cite{vanleeuwen:thesis}(3.3.17(ii)) and \ref{definition:openclosed}, $S_{i,j+1} \in \widetilde{\mC}$.
\begin{itemize}
\item[(a)]  Suppose $S_{i-1,j+1}$ in not in $\OmC$.  Then $r= label
  (S_{i-1,j+1}) < k$. Otherwise $S_{i-1,j}$ and $S_{ij}$ would both
  belong to the same cluster by \cite{vanleeuwen:thesis}(3.3.17(ii)). Since  $S_{i-1,j}$ and $S_{i-1,j+1}$ are in the same cluster by
\cite{vanleeuwen:thesis}(3.3.17(i)) or (3.3.17(iii)), this
  contradicts our assumption. Now \cite{garfinkle1}(1.5.26) forces
$MT(D,T) = \{(k,S_{ij}),(k,S_{i,j+1}) \},$ and since $S_{ij}$ and
$S_{i,j+1}$ both belong to $\widetilde{\mC}$, so must $MT(D,T).$

\item[(b)]  Suppose now that $S_{i-1,j+1} \in \widetilde{\mC}$.  Then the square $S_{i-1,j} \in \OmC$ as well since by \cite{vanleeuwen:thesis}(3.3.17(i)) or (3.3.17(iii)), they both belong to the same cluster.  Now \cite{garfinkle1}(1.5.26) implies $MT(D,T) \subset \{S_{ij}, S_{i-1,j}, S_{i,j+1} \}.$ As all of these squares lie in $\widetilde{\mC}$, we must also have $MT(D,T) \subset \widetilde{\mC}$.
\end{itemize}

\noindent
Case $(ii)$. Suppose $D=\{(k,S_{ij}),(k,S_{i,j+1})\}$ and that the
square $S_{i,j+1}$ is fixed. By
\cite{vanleeuwen:thesis}(3.3.17(ii)) and
\ref{definition:openclosed},  $S_{i,j+2} \in \widetilde{\mC}$.

\begin{itemize}
\item[(a)] Suppose $S_{i-1,j+1}$ is not in $\OmC$.  Then
  $S_{i-1,j+2}$ lies in $|T|$ but not in $\OmC$, as by \cite{vanleeuwen:thesis}(3.3.17(i)) or (3.3.17(iii))  they both belong to the same cluster.  The definition of a cluster forces $r= label (S_{i-1,j+2}) <k$ and \cite{garfinkle1}(1.5.26(ii)) implies $MT(D,T) = \{ S_{i,j+1},S_{i,j+2} \}.$ Since the squares $S_{i,j+1}$ as well as $S_{i,j+2}$ are both contained in $\widetilde{\mC}$, so is $MT(D,T).$
\item[(b)] Suppose $S_{i-1,j+1}$ lies in  $\OmC$.  Then because the
domino $MT(D,T)$ must be a subset of $\{ S_{i,j+1},S_{i,j+2},
S_{i-1,j+1} \}$, it must also be a subset of $\OmC$.
\end{itemize}

\noindent
Case $(iii)$. Suppose $D=\{(k,S_{ij}),(k,S_{i,j+1})\}$ and that the
square $S_{ij}$ is fixed.  Then $S_{i,j-1} \in \OmC$ by
\cite{vanleeuwen:thesis}(3.3.17(iii)).

\begin{itemize}
\item[(a)] Suppose first that $S_{i+1,j-1}$ is not in $\OmC.$  Then $r=
label(S_{i+1,j-1}) >k$ by either
\cite{vanleeuwen:thesis}(3.3.17(ii))   or (3.3.17(iii)). But
\cite{garfinkle1}(1.5.26(iii)) forces $MT(D,T)$ to be precisely
$\{S_{ij},S_{i,j-1} \}$ which is a subset of  $\OmC$.

\item[(b)] If $S_{i+1,j-1} \in \OmC,$ then $S_{i+1,j} \in \widetilde{\mC}$ as
  well, since by \cite{vanleeuwen:thesis}(3.3.17(i)) or (3.3.17(iii)) they either must  belong to the same cluster or $S_{i+1,j}$ is an empty hole or corner.  But by \cite{garfinkle1}(1.5.26(iii)(iv)), $MT(D,T)$ is a subset of $\{S_{ij},S_{i+1,j},S_{i,j-1} \}$, all of whose squares lie in $\widetilde{\mC}$.
\end{itemize}

\noindent
Case $(iv)$.  Suppose $D=\{(k,S_{ij}),(k,S_{i+1,j})\}$ and that the
domino $D$ is of type $I^-$.  The square $S_{i+1,j}$ is then fixed
and $S_{i+1,j-1} \in \OmC$.

\begin{itemize}
\item[(a)] Assume that $S_{i+2,j-1} \in \OmC$.  Then $S_{i+2,j} \in \widetilde{\mC}$.  Since $MT(D,T)$ is the domino $\{S_{i+1,j},S_{i+1,j-1}\}$ or
  $\{S_{i+1,j},S_{i+2,j}\}$. Hence $MT(D,T) \in \OmC$ as both possibilities
  are contained in $\OmC$.

\item[(b)]  Assume $S_{i+2,j-1}$ is not in  $ \OmC$.  We have $r=label(S_{i+2,j-1})>k$, for otherwise $D(r,T)$ and hence $S_{i+2,j-1}$ would lie in  $\OmC$.  But then $MT(D,T) = \{ S_{i+1,j},S_{i+1,j-1} \}$, so it is contained in $\OmC$.
\end{itemize}
These cases describe all possibilities by \ref{theorem:periphery}.
\end{proof}

What remains is to see that the initial cycle $\YC$ is contained
within the cluster $\mC$ itself.  It is enough to show that its
intersection with any closed cluster nested in $\mC$ is empty, as
open clusters cannot be nested. Our proof relies on the notion of
$X$-boxing \cite{garfinkle1}(1.5.2).  We restate the relevant
result.

\begin{proposition}[\cite{garfinkle1}(1.5.9) and (1.5.22)] \label{proposition:xboxed}
Suppose that the dominos $D(k,T)$ and $D(k',T)$ both belong to the
same $X$-cycle.  Then
\begin{enumerate}
        \item   $D(k,T)$ is $X$-boxed iff $MT(D(k,T),T)$ is not $X$-boxed. \label{xboxedi}
        \item  $D(k,T)$ and $D(k',T)$ are both  simultaneously $X$-boxed or not $X$-boxed. \label{xboxedii}
\end{enumerate}
\end{proposition}

\begin{lemma}
If $\mC' \subset \OmC$ is a closed cluster nested in $\mC$, then
$\YC \cap \mC' = \emptyset.$
\end{lemma}

\begin{proof}
It is enough to show that $periphery (\mC') \cap \YC =
\emptyset$, as this forces $\mC' \cap \YC = \emptyset$.  We divide
the problem into a few cases.

Case $(i)$. Suppose $\{type\;\YC, type \;\YCC \} = \{C,D'\}.$  We
  investigate the intersection of  $periphery(\mC')$ with $\YC$.  It cannot contain dominos of types $I^+$ and $I^-$;  because the boxing property
  is constant on cycles according to Proposition \ref{proposition:xboxed}(ii),
 such dominos would have  to be simultaneously $C$ and $D$-boxed,
  which is impossible.  If $D(k,T) \in periphery (\mC') \cap \YCC$ is
  of type $(N)$, $D(k,T)$ and $MT(D(k,T),T)$ are both $C$ and
  $D$'-boxed.  This contradicts Proposition \ref{proposition:xboxed}(i), forcing $periphery (\mC') \cap \YC = \emptyset$.  The proof is identical when the set $\{type\;\YC, type \;\YCC \}$  equals $\{B,D\}$ instead.

Case $(ii)$. Suppose $\{type\;\YC, type \;\YCC \} = \{C,D\}.$  The
  proof is similar to the first case, except this time, dominos of
  type $N$ cannot be simultaneously $C$ and $D$-boxed.  Again, the
  proof is identical when the set $\{type\;\YC, type \;\YCC \}$ equals $\{B,D'\}$ instead.

Case $(iii)$. Suppose $\{type\;\YC, type \;\YCC \} \subset \{B,C\}$ or
  $\{D, D'\}$.  Then by by the definition of cycles, $\YC \cap \YCC =
  \emptyset$.  We know $periphery (\mC') \subset \YCC \subset \OmC'$ by
  Lemma \ref{lemma:cycleinclusterii}, implying again that $periphery
  (\mC') \cap \YC = \emptyset.$
\end{proof}
\section{The $\tau$-Invariant for Orbital Varieties}

A natural question is whether our method of describing orbital
varieties by standard tableaux gives the same
parametrization as \cite{mcgovern:ssmap}. More precisely, if
$\pi: \Irr(\mathcal{F}_u)/A_u \rightarrow \Irr(\mathcal{O}_u \cap
\mathfrak{n})$ is the bijection of \cite{spaltenstein:fixed}, does
the same tableau parameterize both $\mathcal{C} \in
\Irr(\mathcal{F}_u)/A_u$   and its image $\mathcal{V}=\pi(\mathcal{C})$?
Write $\mathcal{T(C)}$ for the domino tableau corresponding to the $A_u$-orbit $\mathcal{C} \in \Irr(\mathcal{F}_u)/A_u$ via the map of the previous section and $\mathcal{T(V)}$ for the domino tableau used to parameterize $\mathcal{V}$ in \cite{mcgovern:ssmap}.

Let $\Pi$ be the set of simple roots in $\mathfrak{g}$.
The $\tau$ invariant, a subset of $\Pi$, is defined for orbital varieties in \cite{joseph:variety} and for components of the Springer fiber in \cite{spaltenstein:book}.  It is constant on each $A_u$ orbit.  For a standard domino tableau $T$, it can be defined in terms of the relative positions of the dominos. We say that a domino $D$ lies {\it higher} than $D'$ in a tableau $T$ iff the rows containing squares of $D$ have indices strictly smaller than the indices of the rows containing squares of $D'$. Then $\tau(T)$ consists of precisely the simple roots $\alpha_i$ whose indices satisfy:

\begin{enumerate}
\item $i=1$ and the domino $D(1,T)$ is vertical,
\item $i>1$ and $D(i-1,T)$ lies higher than $D(i,T)$ in $T$.
\end{enumerate}

According to \cite{garfinkle3}, there is a unique tableau of a given shape within each equivalence class of tableaux generated by the {\it generalized} $\tau$-invariant.  We show

\begin{theorem}
Suppose that $\mathcal{C} \in \Irr(\mathcal{F}_u)/A_u$ and that $\mathcal{V}= \pi(\mathcal{C}).$ Then $$\tau(\mathcal{T(C)}) = \tau(\mathcal{T(V)}).$$
\end{theorem}

\begin{proof}
In fact, we show that all of  following sets are equal.
$$\tau(\mathcal{T(V)})=\tau(\mathcal{V})=\tau(\mathcal{C})= \tau(\mathcal{T(C)}).$$

The first equality follows from \cite{mcgovern:ssmap} and \cite{joseph:variety}.  The second from the definition of $\pi$.  We verify the third.  

Recall the map ${\Phi}: SDT_{op,cl} \rightarrow SDT$ defined in the previous section. We prove that if $\widetilde{T} \in SDT_{op,cl}$ parameterizes the irreducible component $\mathcal{C} \in \Irr\Fu$ in \cite{vanleeuwen:thesis}, then its $\tau$-invariant $\tau(\mathcal{C})$ is precisely the $\tau$-invariant of the standard domino tableau $\Phi(\widetilde{T}) = \mathcal{T(C)}$ as defined above.
The content of the proof is a description of the effect of $\Phi$ on the characterization of the $\tau$-invariant of the components of the Springer fiber given in \cite{spaltenstein:book}.  

That $\alpha_1 \in \tau({\mathcal{C}})$ iff $\alpha_1 \in \tau(\Phi(\widetilde{T})$ is clear in types $B$ and $C$ since $D(1,T)$ never lies within a closed cluster and hence remains unaltered by $\Phi$. In type $D$, the conditions for $\alpha_i,$ when $i \leq 2$, to lie in $\tau(\mathcal{C})$ described by Spaltenstein translate exactly to our conditions for $\alpha_i$ to lie in $\tau(\Phi(\widetilde{T})$.

For $i>1$, suppose that either $D(i,T)$ or $D(i-1,T)$ lies in some $\mathcal{K} \in CC^+(T)$.  If $\mathcal{K}$ contains more than two dominos, then \cite{garfinkle3}(III.1.4) implies that $\alpha_i \in \tau(\mathcal{C}) $ iff $\alpha_i \in \tau(\Phi(\widetilde{T})).$

So suppose that  $\mathcal{K}$ contains exactly two dominos.  If, in fact, $\mathcal{K} = \{D(i), D(i-1)\}$, the simple root $\alpha_i$ must lie in $\tau(\mathcal{C})$. But $D(i-1)$ is higher than $D(i)$ in $MT(\calC, T)$ , implying by the definition of $\Phi$ that $\alpha_i \in \tau(\Phi(\widetilde{T}))$ as well.  The remaining possibility is that only one of the dominos $D(i)$ and $D(i-1)$ lies in the two-domino cluster $\mathcal{K}$. Then the fact that $\alpha_i \in \tau(\mathcal{C})$ iff $\alpha_i \in (\widetilde{T})$ follows by inspection. 
\end{proof}

\providecommand{\bysame}{\leavevmode\hbox to3em{\hrulefill}\thinspace}
\providecommand{\MR}{\relax\ifhmode\unskip\space\fi MR }
\providecommand{\MRhref}[2]{%
  \href{http://www.ams.org/mathscinet-getitem?mr=#1}{#2}
}
\providecommand{\href}[2]{#2}

\end{document}